%% file: Versal.tex
\documentclass[12pt]{smfart}
\usepackage{graphics}

\newtheorem{thm}{Theorem}%[section]
\newtheorem{cor}[thm]{Corollary}

\newtheorem{lem}[thm]{Lemma}

\newtheorem*{theo}{Theorem}
\newtheorem*{pb}{Problem}

\author[Frank Loray]{Frank LORAY}
\address{Frank LORAY (Charg\'e de Recherches au CNRS)\\
IRMAR, UFR de Math\'ematiques,\\
Universit\'e de Rennes 1, Campus de Beaulieu,\\ 
35042 Rennes Cedex (France)}
\email{frank.loray@univ-rennes1.fr}
\urladdr{http://name.math.univ-rennes1.fr/frank.loray}

\title{Versal deformation of the analytic saddle-node}

\begin{document}
\frontmatter
\begin{abstract}
In the continuation of \cite{Prep}, we derive simple forms
for saddle-node singular points of analytic foliations 
in the real or complex plane just by gluing foliated 
complex manifolds. 
We give the versal analytic deformation of the simplest model.
We also derive a unique analytic form 
for those saddle-node having a central manifold.
By this way, we recover and generalize results 
earlier proved by J. \'Ecalle by using mould theory
and partially answer to some questions asked by
J. Martinet and J.-P. Ramis at the end of \cite{MaRa1}.
\end{abstract}
\subjclass{32S65}
\keywords{Normal form, Singularity, Foliation}

\dedicatory{To Jean-Pierre Ramis for his $60^{\text{th}}$ birthday}
\maketitle

\mainmatter

\section*{Introduction}\label{S:intro}

Let $X$ be a germ of analytic vector field 
at the origin of $\mathbb C^2$
$$X=f(x,y)\partial_x+g(x,y)\partial_y,\ \ \ f,g\in\mathbb R\{x,y\}\ 
\text{or}\ \mathbb C\{x,y\}$$
having a singularity at $0$: $f(0)=g(0)=0$.
Consider $\mathcal F$ the germ of singular holomorphic foliation induced 
by the complex integral curves of $X$ near $0$. A question going
back to H. Poincar\'e is the following:

\begin{pb}\rm
Find local coordinates in which the vector field defining the foliation 
has the simplest coefficients.
\end{pb}

For instance, if the vector field $X$ has a linear part (in the matrix form)
$$\left(\begin{matrix} a&b\\ c&d\end{matrix}\right)
=(ax+by)\partial_x+(cx+dy)\partial_y$$
having non zero eigenvalues $\lambda_1,\lambda_2\in\mathbb C$
with non real ratio $\lambda_2/\lambda_1\not\in\mathbb R$,
then H. Poincar\'e proved that the vector field $X$
is actually linear in convenient analytic coordinates. 
In this situation, the eigenvalues $\{\lambda_1,\lambda_2\}$
(resp. the ratio $\lambda_2/\lambda_1$) provide a complete
set of invariants for such vector fields (resp. foliations) 
modulo analytic change of coordinates.

In this paper, we consider {\it saddle-node} singularities of foliations, 
i.e. defined by a vector field with an isolated singularity of the form
$$X=y\partial_y+\text{higher order terms}$$
(one of the eigenvalues is zero). In this situation, H. Dulac proved that,
for arbitrary large $N\in\mathbb N$, there exist analytic coordinates 
in which the foliation is defined by a vector field of the form
\begin{equation}\label{E:preDulac}
X=x^{k+1}\partial_x+y\partial_y+\mu x^ky\partial_y+x^{k+N}f(x,y)\partial_y,\ \ \ 
f\in\mathbb C\{x,y\}
\end{equation}
for unique $k\in\mathbb N^*$ and $\mu\in\mathbb C$. In this sense,
the family of vector fields
\begin{equation}\label{E:formal}
x^{k+1}\partial_x+y\partial_y+\mu x^ky\partial_y
\end{equation}
may be thought as formal normal forms for saddle-node foliations.
The complete analytic classification of those singular points
has been given by J. Martinet and J.-P. Ramis in 1982 (see \cite{MaRa1}),
giving rise to infinitely many invariants additional to the formal
ones $(k,\mu)$ above. Let us recall their main result in the case $k=1$ 
(for simplicity of notations).
Given a saddle-node $\mathcal F$ of the form (\ref{E:preDulac}), 
one can associate the so-called {\it Martinet-Ramis' invariants} 
$(\varphi_0,\varphi_\infty)$ where
\begin{equation}\label{E:cocycle}
\left\{\begin{matrix}
\varphi_0(\zeta)&=&e^{2i\pi\mu}\zeta+\sum_{n\ge2}a_n\zeta^n\ \in\text{Diff}(\mathbb C,0)\\
\varphi_\infty(\zeta)&=&\zeta+t\ \in\mathbb C\ \text{(a translation)}\hfill
\end{matrix}\right.
\end{equation}
The map $\mathcal F\mapsto(\varphi_0,\varphi_\infty)$ is analytic and surjective
onto the set of all pairs satisfying (\ref{E:cocycle}). Any two saddle-nodes 
in the form (\ref{E:preDulac}) with $k=1$ are conjugated by a germ of analytic
diffeomorphism tangent to the identity if and only if the corresponding
invariants coincide. Finally, if we allow linear change of coordinates, 
then one has to consider the pairs (\ref{E:cocycle}) up to an action 
of $\mathbb C^*$ in order to derive the complete analytic classification.
Nevertheless, this final moduli space is singular at the formal models
(\ref{E:formal}) due to an extra symmetry
given by the linear vector field $y\partial_y$. 

The classification above suggests that a generic saddle-node foliation
cannot be defined by a polynomial vector field in any analytic coordinates.
A direct application of our recent work \cite{Prep} provides the following

\begin{thm}\label{T:versal}Let $\mathcal F$ be a germ 
of saddle-node foliation at the origin of $\mathbb{R}^2$ 
(resp. of $\mathbb{C}^2$) having formal invariants $k=1$ and $\mu$ arbitrary.
Then, there exist local analytic coordinates
in which $\mathcal F$ is defined by a vector field of the form
\begin{equation}\label{E:ecalle1}
X_f=x^2\partial_x+y\partial_y+xf(y)\partial_y,\ \ \ f\in\mathbb C\{y\}
\end{equation}
where $f'(0)=\mu$.
\end{thm}

This statement is a particular case of a general simple analytic form
independantly announced by A. D. Brjuno and P. M. Elizarov for all
resonant saddles ($\lambda_2/\lambda_1\in\mathbb Q^-$) 
and saddle-nodes in 1983 (see \cite{Br2,El1}). So 
far, only the case $k=1$ and $\mu=0$ has been proved in 1985 by 
J. \'Ecalle at the end of \cite{Ec}, p.535, presented as an application 
of {\it resurgent functions and mould theory}. In 1994, P. M. Elizarov
considered in \cite{El2} the general saddle-node $\mathcal F$ in Dulac's 
prenormal form (\ref{E:preDulac}) and computed the derivative 
of Martinet-Ramis' modular map $\mathcal F\mapsto(\varphi_0,\varphi_\infty)$ 
at all formal models (\ref{E:formal}). For instance, if we restrict
to the family (\ref{E:ecalle1}) given by Theorem \ref{T:versal},
the modular map is tangent at $X_0=x^2\partial_x+y\partial_y$
to the simpler following one
$$
f(y)=\sum_{n\ge0}a_ny^n\ \mapsto\ 
\left\{\begin{matrix}
\varphi_0(\zeta)&=&e^{2i\pi a_1}\zeta+\sum_{n\ge2}a_n\zeta^n\\
\varphi_\infty(\zeta)&=&\zeta+a_0\hfill
\end{matrix}\right.
$$
In particular, the derivative at $X_0$ is bijective. On the other hand,
our proof of Theorem \ref{T:versal} can be achieved extra complex parameters.
In other words,  
any finite dimensional deformation of $X_0$ inside the family (\ref{E:preDulac})
is analytically conjugated (by an analytic diffeomorphism depending analytically
on the parameters) to a deformation inside (\ref{E:ecalle1}).
In this sense, Theorem \ref{T:versal} provides {\it the versal deformation 
of $X_0$}.

\eject

It is important to mention at this step that the form (\ref{E:ecalle1})
is not unique! Of course, we can modify the functional coefficient $f$
by conjugating the vector field with an homothety $y\mapsto c\cdot y$, 
$c\in\mathbb C^*$. But even if we restrict to tangent-to-the-identity
conjugacies, the form (\ref{E:ecalle1}) is perhaps locally unique at $X_0$
($f\equiv0$), but not globally for the following reason.
By construction (see proof of Theorem \ref{T:versal}), 
the form (\ref{E:ecalle1}) is obtained with $f(0)\neq0$, even if the saddle
node has a {\it central manifold} (see below). For instance, the foliation
from $X_0$ also has another form (\ref{E:ecalle1}) with $f(0)\neq0$.

From preliminary form (\ref{E:preDulac}), we see that $\{x=0\}$
is an invariant curve for the vector field. Tangent to the zero
eigendirection, there is also a ``formal invariant curve'' $\{y=\varphi(x)\}$,
$\varphi\in\mathbb R[[x]]$ or $\mathbb C[[x]]$, which is generically
divergent. Another remarquable result of \cite{MaRa1} is that this formal
invariant curve is actually convergent if, and only if, the translation
part of Martinet-Ramis' invariant is trivial: in the case $k=1$, this means
$\varphi_\infty(\zeta)=\zeta$ (and $t=0$). We then say that the saddle-node
has a {\it central manifold}. For instance, saddle-nodes in the form
(\ref{E:ecalle1}) with $f(0)=0$ have the central manifold $\{y=0\}$.
Conversely, a natural question is

\begin{pb}\rm
Given a saddle-node with a central manifold and formal invariant $k=1$,
is it possible to put it into the form (\ref{E:ecalle1}) with $f(0)=0$
(i.e. simultaneously straightening the central manifold onto $\{y=0\}$) ?
\end{pb}

We provide the complete following answer. 

\begin{thm}\label{T:central}Let $\mathcal F$ be a germ 
of saddle-node foliation at the origin of $\mathbb{R}^2$ 
(resp. of $\mathbb{C}^2$) having formal invariant $k=1$.
Assume that $\mathcal F$ has a central manifold.
Then, there exist local analytic coordinates
in which $\mathcal F$ is defined by a vector field of the form
\begin{equation}\label{E:ecalle2}
X_f=x^2\partial_x+y\partial_y+xf(y)\partial_y,\ \ \ \text{with}\ f(0)=0.
\end{equation}
(and $\mu=f'(0)$) if, and only if, we are in one of the following cases
\begin{enumerate}
\item $\mu\in\mathbb C-\mathbb R^-$, 
\item $\mu<0$ and $\varphi_0$ is linearizable up to conjugacy in
$\text{Diff}(\mathbb C,0)$,
\item $\mu=0$ and Martinet-Ramis' invariants 
$(\tilde\varphi_0^i,\tilde\varphi_\infty^i)_i$ of $\varphi_0$
satisfy: all $\tilde\varphi_0^i$ are linear.
\end{enumerate}
When $\mu\not\in\mathbb Q^-_*$, this form is unique up to homothety $y\mapsto c\cdot y$, 
$c\in\mathbb C^*$.
\end{thm}

\eject

In case (2) of Theorem \ref{T:central}, the linearizability of $\varphi_0$
is automatic (see \cite{Br1}) as soon as $\mu$ belongs to the set $\mathcal B$ 
of Brjuno numbers defined as follows
$$\forall\lambda\in\mathbb R,\ \ \ \lambda\in\mathcal B\ \ \ \Leftrightarrow\ \ \ 
\sum_{n\ge0}{\log(q_{n+1})\over q_n}<\infty$$
(where ${p_n\over q_n}$ stand for succesive truncatures of the continued 
fraction of $\vert\lambda\vert$).
Recall that $\mathcal B$ has full Lebesgue measure in $\mathbb R$. 
But when $\mu\in\mathbb R^-_*\setminus\mathcal B$, the linearizability condition
on $\varphi_0$ is very restrictive (see \cite{Yo}). For instance,
when $\mu\in\mathbb Z-\mathbb N$, this condition forces $\varphi_0$ to be the identity.

The case (3) of Theorem \ref{T:central} has been proved by J. \'Ecalle
in \cite{Ec}, p.539. In this case, $\varphi_0(\zeta)=\zeta+\cdots$ 
is tangent to the identity.
The classification of such map up to conjugacy in $\text{Diff}(\mathbb C,0)$ 
has been given independently by J. \'Ecalle, B. Malgrange and S. M. Voronin 
in 1981. Following the presention of J. Martinet and J.-P. Ramis in \cite{MaRa2}, 
to any element 
$$\varphi(\zeta)=\zeta+\zeta^{n+1}+\cdots\in\text{Diff}(\mathbb C,0)$$
tangent to the identity at the order $n\in\mathbb N$, one associates 
a formal invariant $\lambda\in\mathbb C$ (characterizing the class of $\varphi$
up to formal change of coordinate) and a $n$-uple
$$(\tilde\varphi_0^i,\tilde\varphi_\infty^i)_{i=1,\ldots,n}\in
\left(\text{Diff}(\mathbb C,0)\times\text{Diff}(\mathbb C,\infty)\right)^n$$
$$\text{with}\ \ \ 
\left\{\begin{matrix}
\tilde\varphi_0^i(\zeta)&=&e^{2i\pi\lambda}\zeta+\cdots\in\mathbb C\{\zeta\}\\
\tilde\varphi_\infty^i(\zeta)&=&\zeta+\cdots\hfill\in\mathbb C\{{1\over\zeta}\}
\end{matrix}\right.\ \ \ i=1,\ldots,n
$$
The map $\varphi\mapsto(\tilde\varphi_0,\tilde\varphi_\infty)_i$ is well-defined,
analytic and surjective onto the set of all pairs satisfying conditions above. 
Moreover, any two such maps are conjugated by a change of coordinate tangent 
to the identity if, and only if, the corresponding $n$-uples coincide. Finally,
if we allow linear changes of coordinates, then the conjugacy class of
$\varphi$ within $\text{Diff}(\mathbb C,0)$ is characterized by the $n$-uple
$(\tilde\varphi_0,\tilde\varphi_\infty)_i$ up to an explicit action of $\mathbb C^*$
and a cyclic permutation of indices $i=1,\ldots,n$. Again, the condition
of Theorem \ref{T:central} that all $\tilde\varphi_0$ are linear is very restrictive 
(infinite codimension). In fact, the present work as well as
\cite{Prep} started while we were trying to understand this statement of J. \'Ecalle
in a geometric way. It is now important to say few words about the proofs, which are 
perhaps more interesting than the results above.

\eject

Both forms (\ref{E:ecalle1}) and (\ref{E:ecalle2}) are polynomial
in the variable $x$, and therefore extend analytically along a tubular neighborhood
of $x$-axis, including $x=\infty$: the corresponding singular foliation 
$\mathcal F$ is actually defined on some product $\overline{\mathbb C}\times\Delta$
where $\overline{\mathbb C}=\mathbb C\cup\{\infty\}$ denotes the Riemann sphere
and $\Delta$, a disc on which $f(y)$ is convergent.
For instance, in the form (\ref{E:ecalle1}), the foliation $\mathcal F$
is, apart from the singular point $(x,y)=0$, transversal to the projective line
$L=\overline{\mathbb C}\times\{0\}$ and his tangencies with the vertical projection
$(x,y)\mapsto x$ concentrate along the invariant discs $\{x=0\}$ and $\{x=\infty\}$.
Conversely, it turns out that these geometrical properties essentially characterize 
the form (\ref{E:ecalle1}). 

The proof of Theorem \ref{T:versal} simply consists,
given a germ of saddle-node $\mathcal F$, in gluing it with another foliated complex
domain in order to obtain a germ of $\overline{\mathbb C}$-bundle
$M\simeq\overline{\mathbb C}\times(\mathbb C,0)$ together with a foliation $\mathcal F$
transversal to $\overline{\mathbb C}\times\{0\}$ outside the singular point.
After a convenient uniformization, we obtain the expected form. Although
Theorem \ref{T:versal} is already proved at the end of \cite{Prep},
we reprove it directly in order to be self contained.

By the same way, form (\ref{E:ecalle2}) of Theorem \ref{T:central} is characterized
by the fact that the corresponding foliation $\mathcal F$ extends analytically
on some product $\overline{\mathbb C}\times\Delta$, has central manifold
$L=\overline{\mathbb C}\times\{0\}$ and has exactly one other singular
point along this line, namely at $x=\infty$, with eigenratio 
$\lambda_2/\lambda_1=-\mu$. Since $\varphi_0$ is also the {\it holonomy map}
of the central manifold, necessary and sufficient conditions given 
in Theorem \ref{T:central} exactly coincide with conditions imposed 
by the holonomy of the singular point at $x=\infty$: for instance,
when $\mu<0$ irrational, then $\varphi_0$ must be also the holonomy
of a node with irrational eigenratio $-\mu$ and hence linearizable.
The unicity therefore comes from the unicity of the gluing construction
when the analytic type of the given saddle-node is fixed. In the case
$\mu>0$, we simultaneously obtain a simple form for the saddle at $x=\infty$
(negative eigenratio).

\begin{cor}Let $\mathcal F$ be a germ of saddle foliation at the origin 
of $\mathbb{R}^2$ (resp. of $\mathbb{C}^2$) with eigenratio $-\mu<0$. 
Then there exist local analytic coordinates
in which $\mathcal F$ is defined by a vector field of the form
\begin{equation}\label{E:saddle}
X_f=-x\partial_x+\mu (f(y)+x)y\partial_y,\ \ \ \text{with}\ f(0)=1.
\end{equation}
\end{cor}

This latter form is not unique: when $\mu\in\mathcal B$, all $X_f$ are conjugated.

\eject

For saddle-nodes with formal invariants $k=1$ and $\mu<0$ having a central manifold, 
it is possible to give unique forms by modifying our construction.

\begin{thm}\label{T:penche}Let $\mathcal F$ be a germ 
of saddle-node foliation at the origin of $\mathbb{R}^2$ 
(resp. of $\mathbb{C}^2$) having formal invariant $k=1$.
Assume that $\mathcal F$ has a central manifold.
If $\mu\le0$, then let $n\in\mathbb N$ be such that $\mu+n>0$.
Then, there exist local analytic coordinates
in which $\mathcal F$ is defined by a vector field of the form
\begin{equation}\label{E:brjuno}
X_f=x^2\partial_x+y\partial_y+xyf(x^ny)\partial_y,\ \ \ \text{where}\ f(0)=\mu.
\end{equation}
Moreover, this form is unique up to homothety $y\mapsto c\cdot y$, 
$c\in\mathbb C^*$.
\end{thm}

One could think that we obtain, without analysis, the complete analytic classification 
of saddle-nodes with formal invariant $k=1$ having a central manifold. But the proof
of Theorem \ref{T:penche} uses results of \cite{MaRa1} obtained as corollaries of the 
classification itself. 

Similarly, it is tempting to apply those methods
to the {\it cuspidal} case. The story below has been strongly motivated by questions
and encouragements of Jean-Pierre Ramis.
Recall that, under a generic assumption on the quadratic part, 
the foliation $\mathcal F$ defined by a vector field with nilpotent linear part 
may be written in convenient analytic coordinates (see \cite{CeMo})
$$(2y\partial_x+3 x^2\partial_y) + f(x,y)(2x\partial_x+3y\partial_y)\ \ \ f(0,0)=0.$$
We see that the cusp $\{y^2-x^3=0\}$ is invariant. Following \cite{StZo} (see also \cite{Prep}),
one can analytically reduce $f$ to a function of the single variable $x$.
Still, this form is far to be unique. In fact, one can formally reduce 
$f$ to a unique form $f(x)=x+\tilde f(x^3)$, $\tilde f\in\mathbb C[[x]]$, $\tilde f(0)=0$
(see \cite{Lo}) or analytically reduce any finite jet of $f$ to this form.
But, after experimental tests done by M. Canalis-Durand, F. Michel and M. Teysseire
showing experimental divergence of $\tilde f$ with Gevrey estimates (see \cite{CDMiTe}), 
M. Canalis-Durand and R. Sch\"afke finally established the generic divergence in \cite{CDSc}.
In fact, the equivalent formal reduction to $f(x)=x+\tilde f(h)$ where $h=y^2-x^3$
is the hamiltonian variable (see \cite{Lo}) is Gevrey $1$ and the first singularity 
in the Borel plane is computed in \cite{CDSc}. 

\eject

Clearly, there are at least infinitely many formal invariants but a complete analytic classification
of cuspidal singularities is still missing. In \cite{CeMo}, this problem is reduced to the analytic classification
of a pair of an involution and a trivolution in $\text{Diff}(\mathbb C,0)$ up to simultaneous
conjugacy, but this latter equivalence is not simpler that the initial one. 
One possibility could be to find a unique analytic form. For instance, an easy computation
shows that an analytic form
$$(2y\partial_x+3 x^2\partial_y) + (x+\tilde f(y^2-x^3))(2x\partial_x+3y\partial_y)\ \ \ 
\tilde f(0)=0.$$
defines a foliation $\mathcal F$ that extends analytically along the cusp, and actually
on a tubular neighborhood of the corresponding (singular) rational curve $C$ after
desingularization of the elliptic fibration $\{h=\text{constant}\}$. Along $C$,
apart form the cuspidal singular point, there are exactly $2$ other singular
points for $\mathcal F$, namely a node and a saddle with respective eigenratios
$-6$ and $6$. Clearly the first one is non resonant since it has an invariant curve $C$
in the good direction: it is linearizable and the corresponding holonomy is the identity.

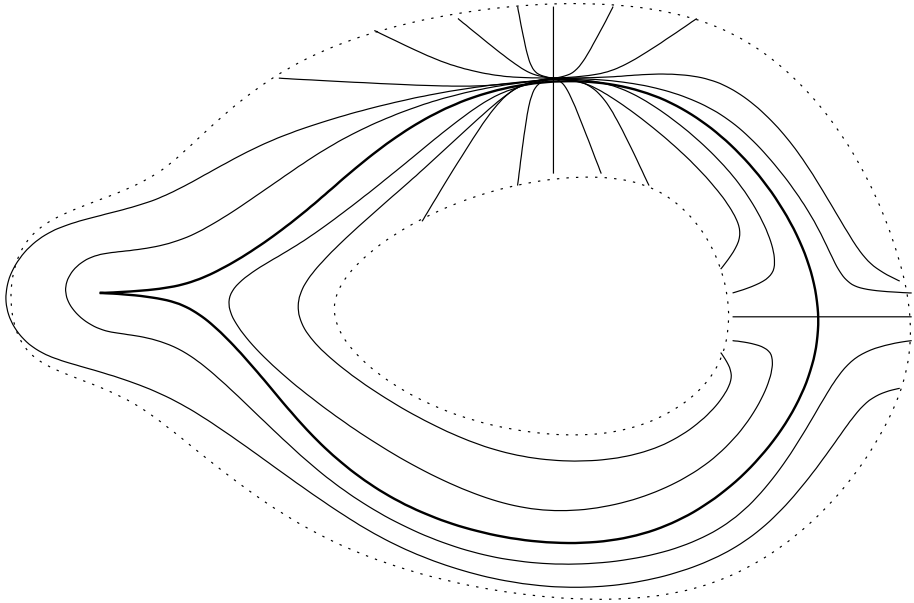
\begin{figure}[htbp]
\begin{center}

\input{Versal1.pstex_t}
 
\caption{The geometry of normal forms from \cite{Lo}}
\label{figure:1}
\end{center}
\end{figure}

On the other hand, the saddle gives no obstruction to the holonomy of the cusp.
Conversely, given a cuspidal foliation $\mathcal F$ (with a generic quadratic part),
one can construct an elliptic fibration having the cusp as a singular fibre 
on which $\mathcal F$ extends with exactly $2$ singular points as above. In order
to do this, we have to construct the elliptic fibration at the same time.
We give in Section \ref{S:gluing} a Lemma which permits to do this (and avoid with
the difficult Savelev's Theorem for the proof of Theorems \ref{T:central} and \ref{T:penche}).
A first difficulty is that the constructed elliptic fibration is not necessarily isotrivial.
But, following Kodaira's classification, the germ of fibration is classified
by an integer, namely the order of modular $\mathcal J$-function at $C$. 
For instance, we have to allow alternate variables like $h={y^2-x^3\over 1+x}$ 
in the non isotrivial case. Now, counting the parameters for instance in the isotrivial case, 
we can write the foliation into the form $X=\partial_t  + f(t,h)h\partial_h$
where $\partial_t$ is a translation along the elliptic fibres, $h\partial_h$ a
multiple of $2x\partial_x+3y\partial_y$, transversal to the regular fibres,
and $f$ an elliptic function on each fibre. By construction, $f$ has $2$ zeroes
and thus $2$ poles: it can be expressed by means of Weierstrass functions
and $4$ coefficients analytic in $h$. After a translation in the fibres,
one can reduce to $3$ analytic coefficients. But at the end, the analytic form
that we could obtain by this way will be not unique because of the resonant node!
Also, we strongly hope that the continuation of \cite{CDSc} will permit 
to provide the good approach for the complete analytic classification of cuspidal
singularities.

\eject

\section{Proof of Theorem \ref{T:versal}}\label{S:versal}

We repeat the geometric construction of \cite{Prep}.
Consider the germ of foliation $\mathcal F_0$ 
defined by a vector field $X_0$ of the form (\ref{E:preDulac})
$$
X_0=x^2\partial_x+y\partial_y+xf(x,y)\partial_y,\ \ \ 
f\in\mathbb C\{x,y\}.
$$
(saddle-node with formal invariant $k=1$ in the form (\ref{E:preDulac}) with $N=0$).
Maybe replacing $y:=x+y$,
the linear part of $X_0$ is given by 
$$\left(\begin{matrix} 0&0\\ c&1\end{matrix}\right)
=(cx+y)\partial_y\ \ \ \text{with}\ \ \ c=f(0)\not=0.$$
Therefore, the vector field $X_0$ is well-defined on the neighborhood 
of any small horizontal disc $\Delta_0=\{\vert x\vert<\varepsilon\}\times\{0\}$, 
$\varepsilon>0$, and transversal to $\Delta_0$ outside the singular point.
Consider inside the horizontal line $L=\overline{\mathbb C}\times\{0\}$
the covering given by $\Delta_0$ and %some complementary disc, say 
$\Delta_\infty=\{\vert x\vert>\varepsilon/2\}\times\{0\}$, and denote
by $C=\Delta_0\cap\Delta_\infty$ the intersection corona.
By the flow-box Theorem,
there exists a unique germ of diffeomorphism of the form
$$\Phi:(\mathbb C^2,C)\to (\mathbb C^2,C)\ ;\ (x,y)\mapsto(\phi(x,y),y),
\ \ \ \phi(x,0)=x$$
straightening $\mathcal F_0$ onto the vertical foliation $\mathcal F_\infty$
(defined by $\partial_y$) at the neighborhood of the corona $C$.
Therefore, after gluing the germs of complex surfaces 
$(\overline{\mathbb C}\times\mathbb C,\Delta_0)$ and 
$(\overline{\mathbb C}\times\mathbb C,\Delta_\infty)$
along the corona by means of $\Phi$, we obtain a germ of smooth complex surface
$S$ along a rational curve $L$ equipped with a singular holomorphic foliation
$\mathcal F$ and a (germ of) rational fibration $y:(S,L)\to(\mathbb C,0)$.
Following \cite{FiGr}, there exists a germ of submersion 
$x:(S,L)\to\overline{\mathbb C}$
completing $y$ into a system of trivializing coordinates: 
$(x,y):(S,L)\to\overline{\mathbb C}\times(\mathbb C,0)$.

\eject

\begin{figure}[htbp]
\begin{center}

\input{Versal2.pstex_t}
 
\caption{Gluing (bi)foliated surfaces}
\label{figure:2}
\end{center}
\end{figure}
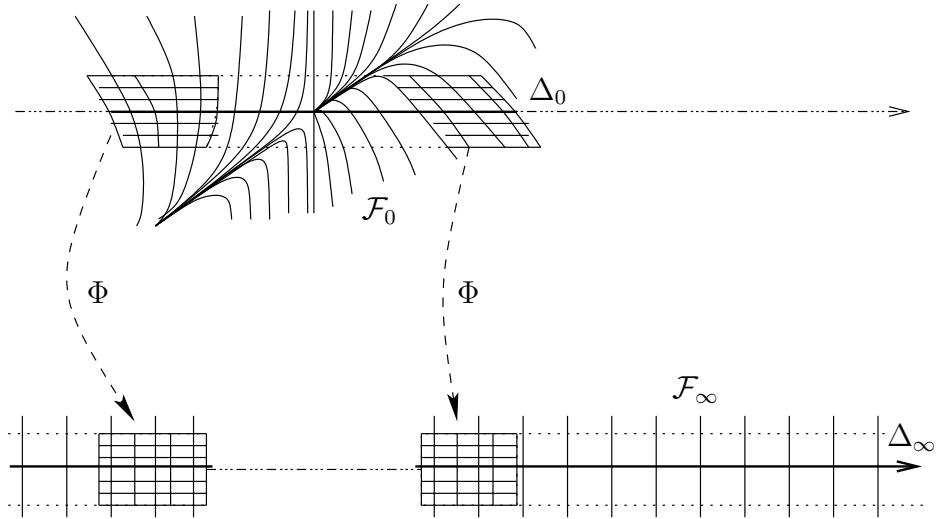

At the neighborhood of any point $p\in L$, the foliation $\mathcal F$ 
is defined by a (non unique) germ of holomorphic vector field,
or equivalently by a unique germ of meromorphic vector field of the form
$$X=f(x,y)\partial_x+\partial_y$$
with $f$ meromorphic at $p$. By unicity, this meromorphic vector field 
is actually globally defined on the neighborhood of $L$ and is therefore
rational in $x$, i.e. $f$ is the quotient of two Weierstrass polynomials. 
For $y$ fixed (close to $0$), the horizontal component $f(x,y)\partial_x$
defines a meromorphic vector field on the corresponding horizontal line  
$\overline{\mathbb C}\times\{y\}$ whose zeroes and poles coincide 
with the tangencies between $\mathcal F$ and the respective vertical 
and horizontal fibrations. By construction, we control the number
of poles: in the second chart, $\mathcal F=\mathcal F_\infty$ is 
transversal to $y$, although in the first chart, $\mathcal F=\mathcal F_0$
has exactly one simple tangency with any horizontal line. It follows
that, for $y$ fixed, the meromorphic vector field $f(x,y)\partial_x$
has exactly $1$ simple pole and thus $3$ zeroes (counted with multiplicity).

\eject

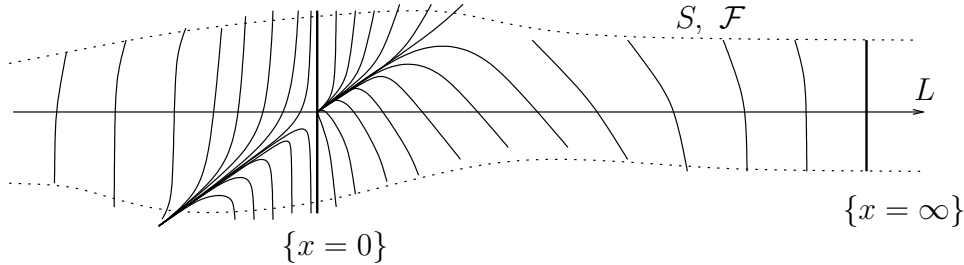
\begin{figure}[htbp]
\begin{center}

\input{Versal3.pstex_t}
 
\caption{Uniformization}
\label{figure:3}
\end{center}
\end{figure}

Of course, in restriction to $L$, the pole vanishes together with one zero
at the singular point of $\mathcal F$. We conclude that the vector field $X$
defining the foliation $\mathcal F$ takes the form 
\begin{equation}\label{E:preparation}
X={f_0(y)+f_1(y)x+f_2(y)x^2+f_3(y)x^3\over g_0(y)+g_1(y)x}\partial_x+\partial_y
\end{equation}
with $f_i,g_j\in\mathbb C\{y\}$. Up to a change of projective horizontal coordinate 
$x:=\{{a(y)x+b(y)\over c(y)x+d(y)}\}$ on $S$, one can assume that $\{x=\infty\}$ 
is a vertical leaf of $\mathcal F$, that $\{x=0\}$ is the invariant curve 
of the saddle-node tangent to the non zero eigendirection and that
$\mathcal F$ has a contact of order $2$ with the vertical fibration 
along $\{x=0\}$ (likely as in the local form (\ref{E:preDulac}) with $k=1$).
Therefore, $f_0,f_1,f_3\equiv0$ and, reminding that $\mathcal F_0$ is a saddle-node
singular point with invariant $k=1$ and $0$-eigendirection transversal to $L$, we also have $f_2(0)\not=0$, $g_1(0)\not=0$,
$g_0(0)=0$ and $g_0'(0)\not=0$. After division, 
$\mathcal F$ is actually defined by a vector field of the form
$$\tilde X=x^2\partial_x+(f(y)x+y g(y))\partial_y,\ \ \ f(0),g(0)\not=0.$$
After change of $y$-coordinate, one may normalize the holomorphic vector field
$y g(y)\partial_y$ to $g(0) y\partial_y$; after division by $g(0)$ and 
linear change of the $x$-coordinate, we finally obtain the form (\ref{E:ecalle1}).

\eject

\section{Proof of Theorem \ref{T:central}}\label{S:central}

Given a saddle-node foliation $\mathcal F$ of the form (\ref{E:ecalle2}),
it is easy to verify that its analytic continuation at the neighborhood 
of the horizontal line $L=\overline{\mathbb C}\times\{0\}$ satisfies
\begin{enumerate}
\item the line $L$ is a global invariant curve for $\mathcal F$, the union
of a smooth leaf together with $2$ singular points;
\item the point $x=0$ is a saddle-node singular point with formal invariants
$k=1$ and $\mu=f'(0)$ and invariant curve $\{xy=0\}$; in particular, 
the sadle-node has an analytic central manifold which is contained in $L$;
\item the point $x=\infty$ is a singular point with eigenratio $-\mu$
and invariant curve $\{x=\infty\}\cup\{y=0\}$;
\item the contact between $\mathcal F$ and the vertical fibration
is double along the invariant curve $\{x=0\}$.
\end{enumerate}

\begin{figure}[htbp]
\begin{center}

\input{Versal4.pstex_t}
 
\caption{Geometry of the second normal form}
\label{figure:4}
\end{center}
\end{figure}
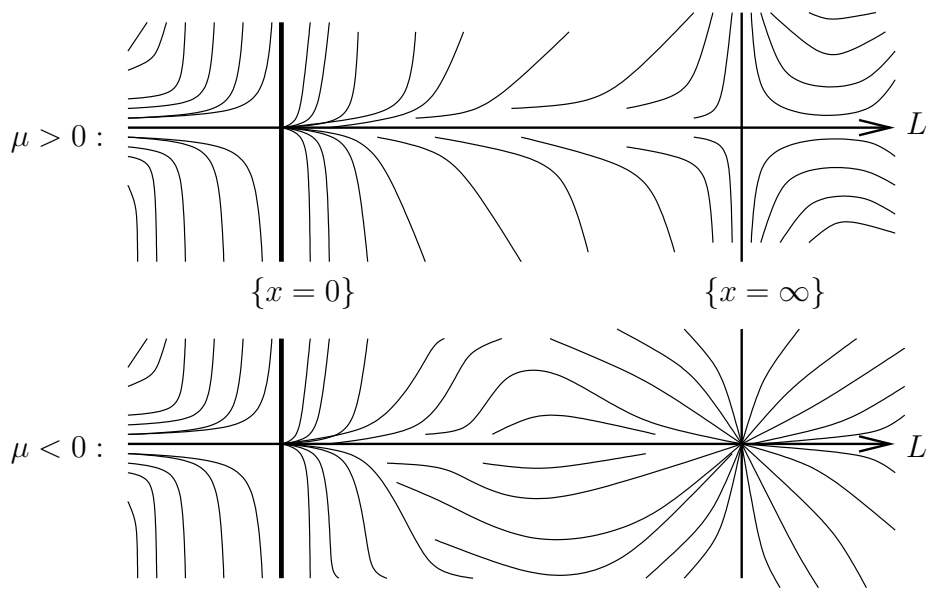

\eject

Conversely, a germ of foliation $\mathcal F$ on 
$\overline{\mathbb C}\times(\mathbb C,0)$ satisfying conditions
above is actually defined by a vector field of the form 
$$\tilde X=x^2\partial_x+(f(y)x+g(y))y\partial_y,\ \ \ f(0),g(0)\not=0.$$
After change of $y$-coordinate, one may normalize the holomorphic vector field
$y g(y)\partial_y$ to $g(0) y\partial_y$; after division by $g(0)$ and 
linear change of the $x$-coordinate, we finally obtain the form (\ref{E:ecalle2}).

\begin{figure}[htbp]
\begin{center}

\input{Versal5.pstex_t}
 
\caption{Holonomy compatibility}
\label{figure:5}
\end{center}
\end{figure}
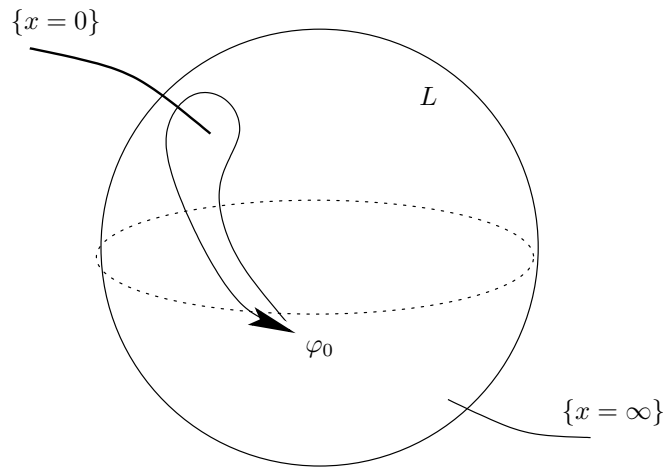

A necessary condition for a saddle-node to admit a form (\ref{E:ecalle2})
is that the holonomy of the central manifold, which actually coincide
with Martinet-Ramis' invariant $\varphi_0$, is also the anti-holonomy
of the invariant curve $L$ around the singular point $x=\infty$.
In the case $\mu<0$, the other singular point is linearizable by Poincar\'e's Theorem 
implying the linearizability of the holonomy map $\varphi_0$.
Here, we use property 3 above and the fact that, in the resonant (non linearizable) case, 
the node has only one irreducible germ of invariant curve.
In the case $\mu=0$, the holonomy $\varphi_0$ is tangent to the identity
and its inverse $\varphi_0^{-1}$ must be the holonomy of the {\it strong manifold}
(the invariant curve tangent to the non zero eigendirection) of a saddle-node
having a central manifold. Following \cite{MaRa1}, this is equivalent to 
condition (3) of Theorem \ref{T:central}.

\eject

We now prove that conditions (1), (2) and (3) of Theorem \ref{T:central}
are sufficient. Like in Section \ref{S:versal}, we start with a germ 
of saddle-node $\mathcal F_0$ defined on the neighborhood of some disc
$\Delta_0$ and glue it with a germ of foliation $\mathcal F_\infty$
along a complementary disc $\Delta_\infty$ in other to obtain a germ
of $2$-dimensional neighborhood $S$ along a rational curve $L$
equipped with a singular foliation $\mathcal F$. The difference 
with Section \ref{S:versal} is that we now glue $\mathcal F_0$ and
$\mathcal F_\infty$ along a common invariant curve, in such a way
that $L$ becomes a global invariant curve for the foliation $\mathcal F$.
In this case, it becomes very technical to do such gluing while preserving
the horizontal foliation. We prefer to glue the foliations with respect 
to the vertical fibration and recover latter the triviality of the neighborhood
by Savelev's Theorem. 

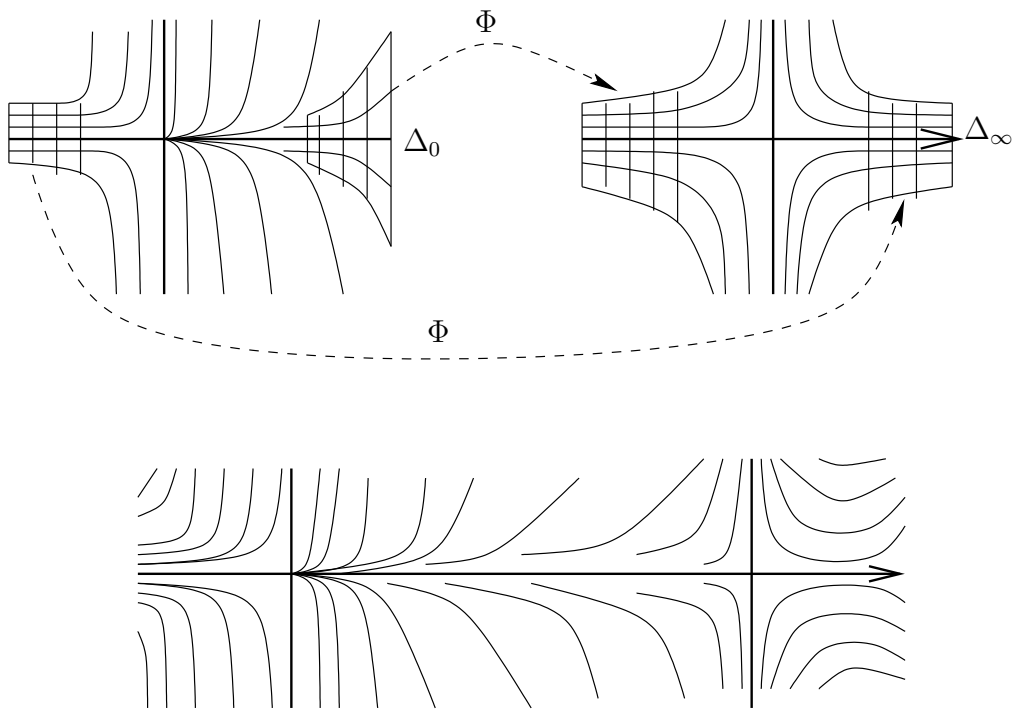
\begin{figure}[htbp]
\begin{center}

\input{Versal6.pstex_t}
 
\caption{Gluing picture}
\label{figure:6}
\end{center}
\end{figure}

\eject

We start with $\mathcal F_0$ into the form
$$
X_0=x^2\partial_x+y\partial_y+xyf(x,y)\partial_y,\ \ \ 
f\in\mathbb C\{x,y\}.
$$
Consider, in local coordinates $(\tilde x=1/x,y)$ at infinity, a germ
of singular foliation $\mathcal F_\infty$ defined by 
$$X_\infty=\tilde x\partial_{\tilde x}-\mu (1+g(\tilde x,y))y\partial_y,\ \ \ 
g(0)=0.$$
Up to a linear change of $\tilde x$-coordinate, one may assume
that $\mathcal F_\infty$ is actually defined on the neighborhood of
$\Delta_\infty$. Obviously, there exists a germ of diffeomorphism 
of the form
$$\Phi:(\mathbb C^2,C)\to (\mathbb C^2,C)\ ;\ (x,y)\mapsto(x,\phi(x,y)),
\ \ \ \phi(x,0)=0$$
straightening $\mathcal F_0$ onto $\mathcal F_\infty$ if, and only if,
the respective holonomy maps around the corona $C=\Delta_0\cap\Delta_\infty$
are conjugated in $\text{Diff}(\mathbb C,0)$.
When $\mu\not\in\mathbb R$, the holonomy map $\varphi_0$ of $\mathcal F_0$
around $C$ (or $x=0$) is hyperbolic and hence linearizable by K\oe nigs' Theorem.
It is therefore enough to choose $X_\infty$ linear.
When $\mu>0$, then the holonomy map $\varphi_0$ can be realized as the holonomy
of a saddle $\mathcal F_\infty$ like above following \cite{MaRa2,PMYo}.
When $\mu<0$ and $\varphi_0$ is linearizable, we obviously realize it with
$X_\infty$ linear. Finally, when $\mu=0$, condition (3) of Theorem \ref{T:central}
is exactly the one to realize $\varphi_0^{-1}$ as the holonomy of a saddle node 
$\mathcal F_\infty$.
After gluing $\mathcal F_0$ and $\mathcal F_\infty$ along $C$, we obtain
a germ of surface $S$ containing a rational curve $L$ which, by Camacho-Sad's
Formula (see \cite{CaSa}), has $0$ self-intersection in $S$. Following
Savelev's Theorem \cite{Sa}, there exists a system of trivializing coordinates: 
$(x,y):(S,L)\to\overline{\mathbb C}\times(\mathbb C,0)$.
Up to a change of trivializing coordinates 
$x:=\{{a(y)x+b(y)\over c(y)x+d(y)}\}$ and $y=\varphi(y)$ on $S$,
one may assume properties (1), (2), (3) and (4) of the begining
of the section all satisfied. Therefore, $\mathcal F$ is defined
by a vector field of the form (\ref{E:ecalle2}).

\eject

Finally, we prove the unicity of form (\ref{E:ecalle2}) in case
$\mu$ is not rational negative. Assume that $\mathcal F$ and 
$\tilde{\mathcal F}$ are of the form (\ref{E:ecalle2}) and 
are analytically conjugated on a neighborhood of $(x,y)=0$.
Following \cite{MaRa1}, they are also conjugated by a germ
of diffeomorphism of the form
$$\Phi_0:(\mathbb C^2,0)\to (\mathbb C^2,0)\ ;\ (x,y)\mapsto(x,\phi_0(x,y))$$
which must preserve the central manifold: $\phi_0(x,0)=0$.
One can extend analytically $\Phi_0$ on a neighborhood of $L-\{x=\infty\}$
in the obvious way, by lifting-path-property. 
The problem is whether $\Phi_0$ extends 
at the other singular point $x=\infty$ or not. If it is the case, 
then we obtain a global diffeomorphism $\Phi$ along $L$
conjugating $\mathcal F$ and $\tilde{\mathcal F}$ and,
by Blanchard's Lemma, permutting the horizontal lines:
$\Phi(x,y)=(x,\phi(y))$. Due to the form (\ref{E:ecalle2}), 
$\phi$ has to commute with $y\partial_y$ and must be linear.
Finally, the fact that $\Phi_0$ extends at the singular point 
at infinity is due to J.-F. Mattei and R. Moussu \cite{MaMo} 
in the case $\mu>0$, to M. Berthier, R. Meziani and P. Sad
in the case $\mu=0$ and we now detail the simpler case $\mathcal F_\infty$
is in the Poincar\'e domain. Recall that property (3) implies
that $\mathcal F_\infty$ is non resonant and hence linearizable 
by a local change of coordinates of the form 
$(\tilde x,y)\mapsto(\tilde x,\phi_\infty(\tilde x,y))$. 
Therefore, we can assume that
$\mathcal F$ and $\tilde{\mathcal F}$ are defined by
$$X_\infty=\tilde x\partial_{\tilde x}-\mu y\partial_y,\ \ \ 
\mu\not\ge0.$$
and that $\Phi_0(\tilde x,y)=(\tilde x,\phi(\tilde x,y))$ 
is a self-conjugacy of $\mathcal F_\infty$ at
the neighborhood of the punctured disc $\Delta^*:=\Delta_\infty-\{\tilde x=0\}$.
The question is, when does $\Phi_0$ coincide with a
symetry of $\mathcal F_\infty$
$$\Phi_\infty(\tilde x,y)=(\tilde x,c\cdot y),\ \ \ c\in\mathbb C^*.$$ 
Of course, this is the case if, and only if, $\phi(\tilde x,y)$ is linear
in $y$. In fact, for $x$ fixed, $\phi(\tilde x,y)$ commutes with the holonomy
$y\mapsto e^{-2i\pi\mu}y$ of $\mathcal F_\infty$ and is therefore linear 
as soon as $\mu$ is not rational. Conversely, in the case $\mu$ is rational,
it is easy to construct examples of $\mathcal F$ and $\tilde{\mathcal F}$
like above that are not globally conjugated and giving rise to non unique
form (\ref{E:ecalle2}).

\eject

\section{Proof of Theorem \ref{T:penche}}

Let us start by blowing-up a saddle-node of the form (\ref{E:ecalle2})  
$$X_f=x^2\partial_x+y\partial_y+xyf(y)\partial_y,\ \ \ f(0)=\mu_0.$$
Along the exceptional diviseur, we have one saddle with eigenratio $-1$
and a saddle-node, given in the chart $(x,t)$, $y=tx$, by
$$\tilde X_f=x^2\partial_x+t\partial_t+xt(f(xt)-1)\partial_t.$$

\begin{figure}[htbp]
\begin{center}

\input{Versal7.pstex_t}
 
\caption{Blowing-up a saddle-node}
\label{figure:7}
\end{center}
\end{figure}
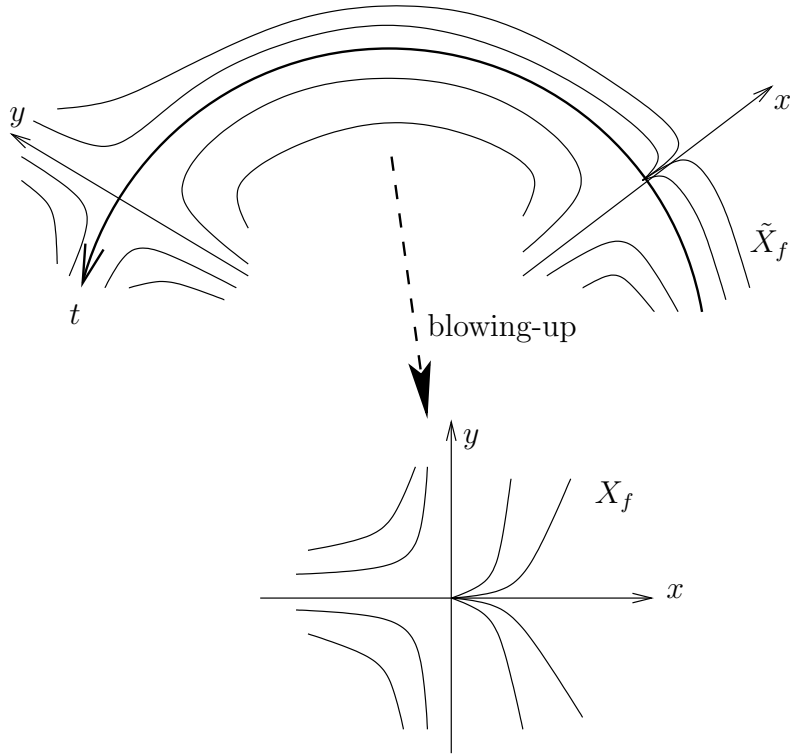

After $n$ successive blowing-up of the saddle-nodes, we obtain
an exceptional divisor like in the picture below where the new
saddle-node takes the form of Theorem \ref{T:penche} with formal
invariant $\mu=\mu_0-n$. All other singular points are saddles
with $-1$ eigenratio. 

\eject

\begin{figure}[htbp]
\begin{center}

\input{Versal8.pstex_t}
 
\caption{After $3$ blowing-up}
\label{figure:8}
\end{center}
\end{figure}
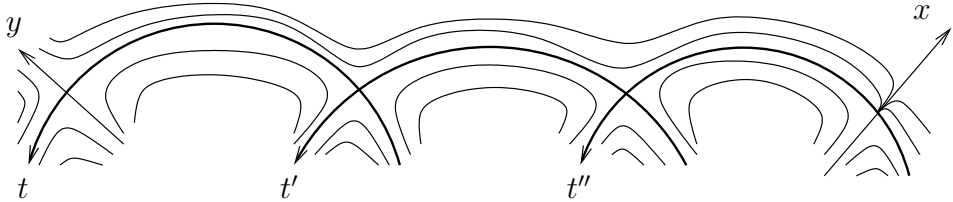

If we are able to prove that any saddle-node $\mathcal F$ 
with invariants $k=1$ and $\mu>-n$ (or $\mu\not\in\mathbb R$) 
having a central manifold is actually the $n$-time blowing-up
of a saddle-node like above, then 
Theorem \ref{T:central} implies Theorem \ref{T:penche}.
To show this, we construct the configuration above by
gluing $\mathcal F$ with successive saddles. Of course,
this is possible because any tangent-to-the-identity 
element of $\text{Diff}(\mathbb C,0)$ can be realized 
as the holonomy map of such saddle (see \cite{MaRa2}).
Then, Grauert's Theorem permits to blow down successively 
the components of the divisor having $-1$ self-intersection.

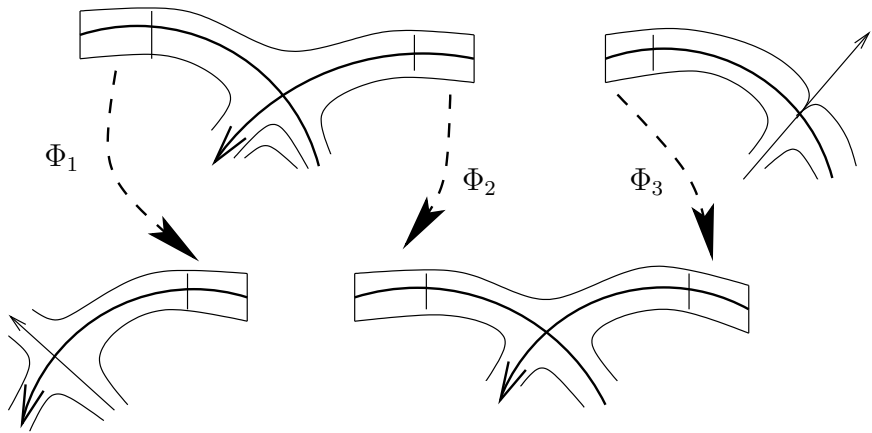
\begin{figure}[htbp]
\begin{center}

\input{Versal9.pstex_t}
 
\caption{Gluing foliations along an exceptional divisor}
\label{figure:9}
\end{center}
\end{figure}

\eject

The unicity of form (\ref{E:brjuno}) also follows from 
that of form (\ref{E:ecalle2}) proved in Section \ref{S:central}.
Indeed, if two such foliations are locally conjugated, 
then the corresponding holonomies of all $-1$ saddle
along the exceptional divisor are also conjugated 
by \cite{MaMo} and, after blowing down, the holonomies
of strong manifold of the sadddle-nodes are so.
By \cite{MaRa1}, the saddle-nodes are also conjugated
after blowing down and they are in the form (\ref{E:ecalle2}).

\section{Gluing Lemmae}\label{S:gluing}

The {\it order of contact} between two germs of regular holomorphic vector fields 
$X_1$ and $X_2$ at $0\in\mathbb C^2$, or between the corresponding
foliations, is by definition the order at $0$ of the determinant $\det(X_1,X_2)$.
For instance, $X_1$ and $X_2$ are transversal if and only if
they have a contact of order $k=0$. Now, if those two foliations share a common
leaf, and if moreover there is no contact between them outside this leaf,
then the contact order $k\in\mathbb N^*$ is constant along this common leaf 
and classifies locally the pair of foliations:

\begin{lem}\label{L:GluingLocal}
Let $\mathcal F$ be a germ of regular analytic foliation at the origin of $\mathbb C^2$
(or $\mathbb R^2$) having the horizontal axis $L_0:\{y=0\}$ as a particular leaf
and having no other contact with the vertical fibration $\{y=\text{constant}\}$:
$\mathcal F$ is defined by a unique function (resp. vector field) of the form
$$F(x,y)=y+y^kxf(x,y)\ \ \ (\text{resp.}\ X=f(x,y)\partial x+y^k\partial_y),\ \ \ 
f(0,0)\not=0$$
where $k\in\mathbb N^*$ denotes the contact order between $\mathcal F$ 
and the fibration.
Then, up to a change of coordinates of the form $\Phi(x,y)=(\phi(x,y),y)$,
the foliation $\mathcal F$ is defined by the function (resp. vector field)
$$F_0(x,y)=y+xy^k\ \ \ (\text{resp.}\ X_0=\partial x+y^k\partial_y).$$
The restriction of $\Phi$ to $L_0$ is the identity if, and only if, $f(0,w)\equiv1$.
Moreover, the normalizing coordinate $\Phi$ is unique once we have decided 
that it fixes the horizontal axis, i.e. $\phi(x,y)=x\tilde\phi(x,y)$.
\end{lem}

\eject

\begin{proof}[Proof]Given $\mathcal F$ as in the statement, 
choose $F(x,y)$ to be the unique function
which is constant on the leaves and has restriction $F(0,y)=y$ on
the horizontal axis: $F(x,y)=y(1+x\tilde F(x,y))$. The assumption
$dF\wedge dy=y^k u(x,y)$, $u(0,0)\not=0$, yields $\tilde
F(x,y)=y^{k-1}f(x,y)$ with $f(0,0)\not=0$, whence the
form $F(x,y)=y+xy^kf(x,y)$. Now, we have
$$F=F_0\circ\Phi_0\ \ \ \text{with}\ \ \ \Phi_0(x,y)=(xf(x,y),y).$$
Thus, $\Phi_0$ is the unique change of $x$-coordinate which conjugates the functions
$F$ and $F_0$; in particular, it conjugates the induced foliations.

Conversely, assume that $\Phi(x,y)=(\phi(x,y),y)$ is conjugating the
foliations respectively induced by $F$ and $F_0$: we have
$$F_0\circ\Phi(x,y)=\varphi\circ F(x,y)\ \ \ \text{with}\ \ \ \varphi(y)=y+y^k\phi(0,y)$$
(the germ $\varphi$ is determined  by the equality restricted to $\{w=0\}$).
If we decompose $f(x,y)=u(x)+y v(x,y)$, we notice that
$\varphi\circ F(x,y)=y+xy^k(u(x)+y\tilde v(x,y))$,
so that $\phi(x,0)=xu(x)=xf(x,0)$.
Finally, if $\phi(x,y)=x\tilde\phi(x,y)$, then $\varphi(y)=y$ and $\Phi$ actually conjugates
the functions: we must have $\Phi=\Phi_0$ whence the unicity.

Now, if $\mathcal F$ is defined by $X=f(x,y)\partial_x+g(x,y)\partial_y$,
assumption gives $dy(X)=g(x,y)=y^k\tilde g(x,y)$ with $f(0,0),\tilde g(0,0)\not=0$. 
After dividing
$X$ by $g$, we can write $X=f(x,y)\partial x+y^k\partial_y$. We have already proved
that any two such foliations (in particular those induced by $X$ and $X_0$) are conjugate
by a unique diffeomorphism of the form $\Phi(x,y)=(x\tilde\phi(x,y),y)$.
Now, if $\Phi(x,y)=(\phi(x,y),y)$ conjugates the foliations respectively induced
by $X$ and $X_0$, it actually conjugates these vector fields. 
In restriction to the trajectory $L_0$, we see that $\phi(x,0)$ conjugates
$\tilde X\vert_{L_0}=f(x,0)\partial_x$ to the constant vector field $\partial_x$.
Therefore, $\phi(x,0)=\int_0^x{1\over f(\zeta,0)}d\zeta$ and
$\phi(x,0)\equiv x$ if, and only if, $f(x,0)\equiv1$.
\end{proof}

\eject

For the next statement, denote by $\Omega\subset(\mathbb C\times\{0\})$ 
a connected open domain inside the horizontal axis.

\begin{lem}\label{L:GluingFibration}
Let $\mathcal F$ and $\mathcal F'$ be regular holomorphic foliations
defined at the neighborhood of $\Omega$ in $\mathbb C^2$
both having $\Omega$ as a particular leaf.
Assume that the contact between each foliation with the horizontal fibration
$\{y=\text{constant}\}$ reduces to $\Omega$, with same order $k\in\mathbb N^*$.
In other words, $\mathcal F$ and $\mathcal F'$ are respectively defined by
vector fields
$$X=f(x,y)\partial_x+y^k\partial_y\ \ \ \text{and}\ \ \ X'=f'(x,y)\partial_x+y^k\partial_y$$
where $f$ and $f'$ are non vanishing functions in the neighborhood of $\Omega$.
Then, $\mathcal F$ and $\mathcal F'$ are conjugated in a neighborhood of $\Omega$
by a diffeomorphism of the form $\Phi(x,y)=(x+y\phi(x,y),\psi(y))$ (fixing $\Omega$)
if, and only if, the two following conditions hold
\begin{enumerate}
\item $f(x,0)\equiv f'(x,0)$;
\item the respective holonomies $\varphi$ and $\varphi'$ of $\mathcal F$ 
and $\mathcal F'$
along $\Omega$ are analytically conjugated: $\psi\circ\varphi=\varphi'\circ\psi$.
\end{enumerate}
\end{lem}

\begin{proof}[Proof]Following Lemma \ref{L:GluingLocal}, condition (1) is the necessary
and sufficient condition for the existence of local conjugacies $\Phi=(y\phi(x,y),y)$
between $\mathcal F$ and $\mathcal F'$ at the neighborhood of any point $w_0\in\Omega$.
Fix one of these points and consider the respective holonomy maps $\varphi$ and $\varphi'$
computed on the transversal $T:\{x=x_0\}$ in the variable $y$. By condition (2),
up to conjugate, say $\mathcal F'$, by a diffeomorphism of the form $(x,\psi(y))$,
we may assume without loss of generality $\varphi(y)=\varphi'(y)$.
We start with the local diffeomorphism $\Phi(x,y)=(y\phi(x,y),y)$ given by 
Lemma \ref{L:GluingLocal}
conjugating the foliations and fixing $T$. Since $\Phi$ conjugates the corresponding
vector fields $X$ and $X'$, it extends analytically along the whole of $\Omega$ by the formula
$\Phi(p):=\Phi^{-t}_{X'}\circ\Phi\circ\Phi^t_X(p)$. The condition $\varphi(y)=\varphi'(y)$
implies that $\Phi$ is uniform.
\end{proof}

\eject

Here is another gluing Lemma for pairs of regular foliations
$\mathcal F$ and $\mathcal G$ at the neighborhood of the common leaf $\Omega$
(before, $\mathcal G$ was the horizontal fibration).

\begin{lem}Let $\mathcal F$ and $\mathcal G$ (resp. $\mathcal F'$ and $\mathcal G'$)
be regular holomorphic foliations defined at the neighborhood of the corona $\Omega$ in $\mathbb C^2$
both having $\Omega$ as a regular leaf.
Assume that the contact between $\mathcal F$ and $\mathcal G$ (resp. $\mathcal F'$ and $\mathcal G'$)
reduces to $\Omega$, with same order $k\in\mathbb N^*$.
In other words, the foliations above are respectively defined by
vector fields
$$X=\partial_x+yf(x,y)\partial_y\ \ \ \text{and}\ \ \ Y=X+y^kg(x,y)\partial_y,$$
$$(\text{resp.}\ \ \ X'=\partial_x+yf'(x,y)\partial_y\ \ \ \text{and}\ \ \ 
Y'=X'+y^kg'(x,y)\partial_y)$$
where $g$ and $g'$ are non vanishing functions in the neighborhood of $\Omega$.
Then, $\mathcal F$ and $\mathcal G$ are simultaneously conjugated to $\mathcal F'$ and $\mathcal G'$
in a neighborhood of $\Omega$ by a diffeomorphism of the form $\Phi(x,y)=(x+y\phi(x,y),y\psi(x,y))$
(fixing point-wise $\Omega$) if, and only if, the two conditions hold
\begin{enumerate}
\item  for any (and for all) $x_0\in\Omega$, we have
$${g(x,0)\over\exp(-\int_{x_0}^x f(\zeta,0)d\zeta)}\equiv
 {g'(x,0)\over\exp(-\int_{x_0}^x f'(\zeta,0)d\zeta)};$$
\item the respective pairs of holonomies $(\varphi_\mathcal F,\varphi_\mathcal G)$ and
$(\varphi_{\mathcal F'},\varphi_{\mathcal G'})$ along $\Omega$
are simultaneously analytically conjugated:
$\psi\circ\varphi_\mathcal F=\varphi_{\mathcal F'}\circ\psi$ and
$\psi\circ\varphi_\mathcal G=\varphi_{\mathcal G'}\circ\psi$.
\end{enumerate}
\end{lem}

\begin{proof}[Proof]It is similar to that of the previous Lemma.
Up to a change of coordinate $y:=\psi(y)$ (which does not affect neither $f(x,0)$, nor $g(x,0)$
and hence preserves equality (i)), we may assume that holonomies
$(\varphi_\mathcal F,\varphi_\mathcal G)\equiv(\varphi_{\mathcal F'},\varphi_{\mathcal G'})$
actually coincide on a transversal $T:\{x=x_0\}$. We just detail that condition (i)
exactly provides the existence of local conjugacies between the given pairs of foliations
fixing point-wise $\Omega$; the unique conjugacy
fixing $T$ will extend uniformly along $\Omega$ by (ii).

At the neighborhood of any point $(x_0,0)$ of the corona $\Omega$,
say $x_0=0$ for simplicity, we preliminary conjugate $X$ to $X_0=\partial_x$
by respective local changes of $y$-coordinate $\Psi(x,y)=(x,y\psi(x,y))$, $\psi(0,0)\not=0$
$$ \Psi_*X=X_0=\partial_x \ \ \ \text{and}\ \ \
  \Psi_*Y=Y_0=z^kg_0(x,y)\partial_z+\partial_w$$
Doing the same with the pair $X'$ and $Y'$, we see by
Lemma \ref{L:GluingLocal} that the corresponding pairs of foliations are conjugated
by a diffeomorphism fixing point-wise $\{y=0\}$ if, and only if,
the differential form $\omega=g_0(x,0)dx$ along $\Omega$ coincide
with the corresponding one $\omega'=g_0'(x,0)dx$
for $X_0'=X_0$ and $Y_0'=\partial_x+y^kg_0(x,y)\partial_y$.
This $1$-form $\omega$ can be redefined in the following intrinsic way:
the holonomy of $\mathcal G$ between two transversal cross-sections $T_0$
and $T_1$ computed in any coordinate $y$ which is $\mathcal F$-invariant
(here $\mathcal F$ is defined by $\partial_x$) takes the form
$$\varphi(y)=y+\left(\int_{x_0}^{x_1}\omega\right)y^k+(\text{higher order terms})$$
where $(x_i,0):=T_i\cap\Omega$, $i=0,1$. Since
$$\Psi^*X_0=\partial_x-{\psi_x\over\psi+y\psi_y}y\partial_y\ \ \ \text{and}\ \ \
\Psi^*Y_0=\Psi^*X_0+{g_0\over\psi+y\psi_y}y^k\partial_y,$$
($\psi_x$ and $\psi_y$ are partial derivatives of $\psi$) we derive in restriction to $\Omega$
$$f(x,0)=-{\psi_x(x,0)\over\psi(x,0)}\ \ \ \text{and}\ \ \
g_0(x,0)=\psi(x,0)\cdot g(x,0)$$
yielding the formula for the local invariant of our conjugacy problem
$$\omega={g(x,0)\over\exp(-\int_{x_0}^x f(\zeta,0)d\zeta)}dx.$$
\end{proof}

\eject

\section{Complements}\label{S:complements}

Let $X$ be a germ of analytic vector field 
at the origin of $\mathbb C^2$
$$X=f(x,y)\partial_x+g(x,y)\partial_y,\ \ \ f,g\in\mathbb R\{x,y\}\ 
\text{or}\ \mathbb C\{x,y\}$$
having a singularity at $0$, $f(0)=g(0)=0$,
and consider its linear part
$$X_1=(ax+by)\partial_x+(cx+dy)\partial_y
=\left(\begin{matrix} a&b\\ c&d\end{matrix}\right)$$
(written in matrix form) where $X=X_1+\cdots$ (quadratic terms).
Given an analytic diffeomorphism $\Phi=(\phi_1,\phi_2):(\mathbb C^2,0)\to(\mathbb C^2,0)$,
we have
$$\Phi^*X=D\Phi^{-1}\cdot X\circ\Phi 
=\left(\begin{matrix}\phi_{1x}&\phi_{1y} \\ \phi_{2x}&\phi_{2y}\end{matrix}\right)^{-1}
\cdot\left(\begin{matrix}f\circ\Phi \\ g\circ\Phi\end{matrix}\right).$$
In particular, the linear part of $\Phi^*X$ is given by
$$D\Phi(0)^*X_1=
\left(\begin{matrix}\phi_{1x}(0)&\phi_{1y}(0) \\ \phi_{2x}(0)&\phi_{2y}(0)\end{matrix}\right)^{-1}
\cdot\left(\begin{matrix}  a& b\\  c& d\end{matrix}\right)
\cdot \left(\begin{matrix}\phi_{1x}(0)&\phi_{1y}(0) \\ \phi_{2x}(0)&\phi_{2y}(0)\end{matrix}\right)$$
and the eigenvalues $\lambda_1,\lambda_2\in\mathbb C$ of (the linear part of)
the vector field $X$ are invariant under analytic (or formal) change of coordinates.
When one of the eigenvalues is not zero, say $\lambda_1$, 
we denote by $\lambda:=\lambda_2/\lambda_1\in\mathbb C$ the eigenratio.
Since $h\cdot X$ defines the same foliation for any $h\in\mathbb C\{x,y\}$,
$h(0)\not=0$, the eignevalues $\lambda_i$ are no more invariants for the foliation 
but the eigenratio $\lambda$ is. In this situation, a famous result of H. Poincar\'e,
improved in the dimension $2$ case by H. Dulac, yield

\begin{theo}[Poincar\'e-Dulac]
If $\lambda\not\in\mathbb{R}^-$, 
then, up to a (real or complex) analytic change of coordinates,
the vector field becomes either linear diagonal
$$X=\lambda_1x\partial_x+\lambda_2 y\partial_y$$
or, in the resonant node case $\{\lambda_1,\lambda_2\}=\{\lambda,k\lambda\}$,
$k\in\mathbb{N}^*$ and $\lambda\in\mathbb C^*$, possibly of the form
(including the non diagonal linear case for $n=1$)
$$X=\lambda(kx+y^k)\partial_x+\lambda y\partial_y$$
or, in the real hyperbolic focus case $\{\lambda_1,\lambda_2\}=\{a\pm ib\}$, $a,b\in\mathbb R^*$
$$X=(ax-by)\partial_x+(bx+ay)\partial_y.$$
\end{theo}

Those normal forms are unique (up to a permutation of $\lambda_1$ and $\lambda_2$)
so that this statement cannot be improved. In the case $\lambda\in\mathbb{R}^-$, 
analogous simple polynomial models are obtained for the subjacent foliation
after a formal change of coordinates, that is a pair $\Phi=(\phi_1,\phi_2)$
where $\phi_1,\phi_2\in\mathbb C[[x,y]]$ and $\det D\Phi(0)\not=0$.

\begin{theo}[Poincar\'e-Dulac]
If $\lambda\in\mathbb{R}^-$, then, up to a formal change of coordinates 
and multiplication by a non vanishing formal power series,
the vector field becomes either linear
$$X=x\partial_x+\lambda y\partial_y$$
or, in the resonant saddle case $\lambda=-{p\over q}\in\mathbb{Q}^*_-$,
possibly of the form
$$X=p x\partial_x-\left(q+(x^py^q)^k+\mu(x^py^q)^{2k}\right)y\partial_y$$
or, in the saddle-node case $\lambda=0$, of the form
$$X=x\partial_x+\left(y^k+\mu y^{2k}\right)y\partial_y$$
for unique positive integer $k\in\mathbb{N}^*$ and scalar $\mu\in\mathbb{C}$.
\end{theo}

See \cite{Me} or \cite{Te} for formal models of the corresponding vector fields.
When $\lambda$ belongs to the set $\mathcal B\cap\mathbb R^-$ of Brjuno numbers,
A. D. Brjuno proved in 1971 (see \cite{Br1}) that the vector field $X$ is actually 
analytically linearisable. The remaining exceptional
cases split into diophantine saddles 
($\lambda\in\mathbb R^-\setminus(\mathbb Q^-\cup\mathcal B)$)
resonant saddles ($\lambda\in\mathbb Q^-\setminus\{0\}$) and 
saddle-nodes ($\lambda=0$). In the resonant case ($\lambda\in\mathbb R^-$), 
Yu. S. Ilyashenko proved in 1981 (see \cite{Il})
that the formal change of coordinates leading to the Poincar\'e-Dulac normal form
above is divergent as a rule. The complete analytic classification 
of saddle-nodes and resonant saddles has been given by J. Martinet and J.-P. Ramis
in 1982 and 1983 (see \cite{MaRa1,MaRa2}) giving rise to infinitely many invariants 
additional to the formal ones $(-{p\over q},k,\mu)$ above. 
See \cite{GrVo,Gr} and \cite{MeVo,Te} for the corresponding classification 
of vector fields. For diophantine saddles, although a complete classification
is still missing, J.-C. Yoccoz proved in 1988 (see \cite{Yo}) that the moduli space
between analytic and formal classification is also very huge (infinite dimension).

\eject

Let us now consider the following family of saddle-nodes ($\varepsilon>0$)
$$X_\varepsilon={x^2\over 1+\mu x}\partial_x+y\partial_y+\varepsilon f(x,y)y\partial_y,\ \ \ 
f=\sum_{m\le0,\ n\le-1}f_{m,n}x^my^n\in{1\over y}\mathbb C\{x,y\}$$
with formal invariants $k=1$ and $\mu\in\mathbb C$: $f_{0,0}=f_{0,1}=f_{1,1}=0$. 
Consider also the associate Martinet-Ramis' invariants (depending on $\varepsilon$)
$$
\varphi_0(\zeta)=e^{2i\pi\mu}\zeta+\sum_{n>0}\varphi_n\zeta^{n+1}\ \ \ \text{and}\ \ \ 
\varphi_\infty(\zeta)=\zeta+t
$$ 
Then, the main result of \cite{El2} reads

\begin{theo}[Elizarov]The derivative (in the sense of G\^ateau)
of Martinet-Ramis'moduli at $\varepsilon=0$ is given by
$${d\varphi_n\over\partial\varepsilon}\vert_{\varepsilon=0}
=n^{\mu n-1}e^{-2i\pi n\mu}\sum_{m>0}{m\over\Gamma(1+m+\mu n)}f_{m,n}(-n)^m$$
$$\text{and}\ \ \ {dt\over\partial\varepsilon}\vert_{\varepsilon=0}
=(-1)^{-\mu}e^{2i\pi\mu}\sum_{m>0}{m\over\Gamma(1+m-\mu)}f_{m,-1}(-n)^m$$
where $\Gamma$ is the Euler's Gamma Function.
\end{theo}

Finally, the analytic form announced in \cite{Br2} is the following

{\it Given a saddle-node $\mathcal F$ with formal invariants $k=1$ 
and $\mu\in\mathbb C$ and a slope $0<s\le+\infty$, 
one can find analytic coordinates in which $\mathcal F$
is defined by a vector field of the form
$$x^2\partial_x+y\partial_y+
x\left(f_0+\mu y+\sum_{(m,n)\in E_s}f_{m,n}x^my^{n+1}\right)\partial_y$$
provided that the set of exponents $E_s=\{(m,n);n>0,{n\over s}+1\le m<{n\over s}+2\}$ 
does not intersect the resonances $\{(m,n);m+\mu n\in-\mathbb N\}$.}

\eject

\end{document}

%% file: Versal1.pstex_t
\begin{picture}(0,0)%
\includegraphics{Versal1.pstex}%
\end{picture}%
\setlength{\unitlength}{3947sp}%
\begingroup\makeatletter\ifx\SetFigFont\undefined%
\gdef\SetFigFont#1#2#3#4#5{%
  \reset@font\fontsize{#1}{#2pt}%
  \fontfamily{#3}\fontseries{#4}\fontshape{#5}%
  \selectfont}%
\fi\endgroup%
\begin{picture}(5713,3769)(0,-3799)
\end{picture}

%% file: Versal2.pstex_t
\begin{picture}(0,0)%
\includegraphics{Versal2.pstex}%
\end{picture}%
\setlength{\unitlength}{3947sp}%
\begingroup\makeatletter\ifx\SetFigFont\undefined%
\gdef\SetFigFont#1#2#3#4#5{%
  \reset@font\fontsize{#1}{#2pt}%
  \fontfamily{#3}\fontseries{#4}\fontshape{#5}%
  \selectfont}%
\fi\endgroup%
\begin{picture}(5806,4138)(1114,-4412)
\put(3376,-1636){\makebox(0,0)[lb]{\smash{\SetFigFont{12}{14.4}{\familydefault}{\mddefault}{\updefault}$\mathcal F_0$}}}
\put(6676,-3061){\makebox(0,0)[lb]{\smash{\SetFigFont{12}{14.4}{\familydefault}{\mddefault}{\updefault}$\Delta_\infty$}}}
\put(1651,-2161){\makebox(0,0)[lb]{\smash{\SetFigFont{12}{14.4}{\familydefault}{\mddefault}{\updefault}$\Phi$}}}
\put(5326,-2761){\makebox(0,0)[lb]{\smash{\SetFigFont{12}{14.4}{\familydefault}{\mddefault}{\updefault}$\mathcal F_\infty$}}}
\put(3976,-2161){\makebox(0,0)[lb]{\smash{\SetFigFont{12}{14.4}{\familydefault}{\mddefault}{\updefault}$\Phi$}}}
\put(4426,-886){\makebox(0,0)[lb]{\smash{\SetFigFont{12}{14.4}{\familydefault}{\mddefault}{\updefault}$\Delta_0$}}}
\end{picture}

%% file: Versal3.pstex_t
\begin{picture}(0,0)%
\includegraphics{Versal3.pstex}%
\end{picture}%
\setlength{\unitlength}{3947sp}%
\begingroup\makeatletter\ifx\SetFigFont\undefined%
\gdef\SetFigFont#1#2#3#4#5{%
  \reset@font\fontsize{#1}{#2pt}%
  \fontfamily{#3}\fontseries{#4}\fontshape{#5}%
  \selectfont}%
\fi\endgroup%
\begin{picture}(5799,3163)(1114,-3437)
\put(6826,-886){\makebox(0,0)[lb]{\smash{\SetFigFont{12}{14.4}{\familydefault}{\mddefault}{\updefault}$L$}}}
\put(5326,-436){\makebox(0,0)[lb]{\smash{\SetFigFont{12}{14.4}{\familydefault}{\mddefault}{\updefault}$S,\ \mathcal F$}}}
\put(2851,-1861){\makebox(0,0)[lb]{\smash{\SetFigFont{12}{14.4}{\familydefault}{\mddefault}{\updefault}$\{x=0\}$}}}
\put(6376,-1636){\makebox(0,0)[lb]{\smash{\SetFigFont{12}{14.4}{\familydefault}{\mddefault}{\updefault}$\{x=\infty\}$}}}
\end{picture}

%% file: Versal4.pstex_t
\begin{picture}(0,0)%
\includegraphics{Versal4.pstex}%
\end{picture}%
\setlength{\unitlength}{3947sp}%
\begingroup\makeatletter\ifx\SetFigFont\undefined%
\gdef\SetFigFont#1#2#3#4#5{%
  \reset@font\fontsize{#1}{#2pt}%
  \fontfamily{#3}\fontseries{#4}\fontshape{#5}%
  \selectfont}%
\fi\endgroup%
\begin{picture}(5626,3649)(451,-3073)
\put(4801,-1261){\makebox(0,0)[lb]{\smash{\SetFigFont{12}{14.4}{\familydefault}{\mddefault}{\updefault}$\{x=\infty\}$}}}
\put(6076,-211){\makebox(0,0)[lb]{\smash{\SetFigFont{12}{14.4}{\familydefault}{\mddefault}{\updefault}$L$}}}
\put(6076,-2236){\makebox(0,0)[lb]{\smash{\SetFigFont{12}{14.4}{\familydefault}{\mddefault}{\updefault}$L$}}}
\put(1951,-1261){\makebox(0,0)[lb]{\smash{\SetFigFont{12}{14.4}{\familydefault}{\mddefault}{\updefault}$\{x=0\}$}}}
\put(451,-286){\makebox(0,0)[lb]{\smash{\SetFigFont{12}{14.4}{\familydefault}{\mddefault}{\updefault}$\mu>0:$}}}
\put(451,-2236){\makebox(0,0)[lb]{\smash{\SetFigFont{12}{14.4}{\familydefault}{\mddefault}{\updefault}$\mu<0:$}}}
\end{picture}

%% file: Versal5.pstex_t
\begin{picture}(0,0)%
\includegraphics{Versal5.pstex}%
\end{picture}%
\setlength{\unitlength}{3947sp}%
\begingroup\makeatletter\ifx\SetFigFont\undefined%
\gdef\SetFigFont#1#2#3#4#5{%
  \reset@font\fontsize{#1}{#2pt}%
  \fontfamily{#3}\fontseries{#4}\fontshape{#5}%
  \selectfont}%
\fi\endgroup%
\begin{picture}(3653,2878)(376,-2334)
\put(376,420){\makebox(0,0)[lb]{\smash{\SetFigFont{10}{12.0}{\familydefault}{\mddefault}{\updefault}$\{x=0\}$}}}
\put(2226,-1610){\makebox(0,0)[lb]{\smash{\SetFigFont{10}{12.0}{\familydefault}{\mddefault}{\updefault}$\varphi_0$}}}
\put(3837,-2028){\makebox(0,0)[lb]{\smash{\SetFigFont{10}{12.0}{\familydefault}{\mddefault}{\updefault}$\{x=\infty\}$}}}
\put(2942,-58){\makebox(0,0)[lb]{\smash{\SetFigFont{10}{12.0}{\familydefault}{\mddefault}{\updefault}$L$}}}
\end{picture}

%% file: Versal6.pstex_t
\begin{picture}(0,0)%
\includegraphics{Versal6.pstex}%
\end{picture}%
\setlength{\unitlength}{3947sp}%
\begingroup\makeatletter\ifx\SetFigFont\undefined%
\gdef\SetFigFont#1#2#3#4#5{%
  \reset@font\fontsize{#1}{#2pt}%
  \fontfamily{#3}\fontseries{#4}\fontshape{#5}%
  \selectfont}%
\fi\endgroup%
\begin{picture}(6044,4405)(1179,-3734)
\put(3826,-1411){\makebox(0,0)[lb]{\smash{\SetFigFont{12}{14.4}{\familydefault}{\mddefault}{\updefault}$\Phi$}}}
\put(4126,539){\makebox(0,0)[lb]{\smash{\SetFigFont{12}{14.4}{\familydefault}{\mddefault}{\updefault}$\Phi$}}}
\put(3676,-211){\makebox(0,0)[lb]{\smash{\SetFigFont{12}{14.4}{\familydefault}{\mddefault}{\updefault}$\Delta_0$}}}
\put(7201,-136){\makebox(0,0)[lb]{\smash{\SetFigFont{12}{14.4}{\familydefault}{\mddefault}{\updefault}$\Delta_\infty$}}}
\end{picture}

%% file: Versal7.pstex_t
\begin{picture}(0,0)%
\includegraphics{Versal7.pstex}%
\end{picture}%
\setlength{\unitlength}{3947sp}%
\begingroup\makeatletter\ifx\SetFigFont\undefined%
\gdef\SetFigFont#1#2#3#4#5{%
  \reset@font\fontsize{#1}{#2pt}%
  \fontfamily{#3}\fontseries{#4}\fontshape{#5}%
  \selectfont}%
\fi\endgroup%
\begin{picture}(4824,4722)(1564,-4723)
\put(4201,-2086){\makebox(0,0)[lb]{\smash{\SetFigFont{12}{14.4}{\familydefault}{\mddefault}{\updefault}blowing-up}}}
\put(5701,-3736){\makebox(0,0)[lb]{\smash{\SetFigFont{12}{14.4}{\familydefault}{\mddefault}{\updefault}$x$}}}
\put(4426,-2761){\makebox(0,0)[lb]{\smash{\SetFigFont{12}{14.4}{\familydefault}{\mddefault}{\updefault}$y$}}}
\put(6376,-661){\makebox(0,0)[lb]{\smash{\SetFigFont{12}{14.4}{\familydefault}{\mddefault}{\updefault}$x$}}}
\put(1576,-736){\makebox(0,0)[lb]{\smash{\SetFigFont{12}{14.4}{\familydefault}{\mddefault}{\updefault}$y$}}}
\put(5251,-3136){\makebox(0,0)[lb]{\smash{\SetFigFont{12}{14.4}{\familydefault}{\mddefault}{\updefault}$X_f$}}}
\put(6226,-1561){\makebox(0,0)[lb]{\smash{\SetFigFont{12}{14.4}{\familydefault}{\mddefault}{\updefault}$\tilde X_f$}}}
\put(1951,-2011){\makebox(0,0)[lb]{\smash{\SetFigFont{12}{14.4}{\familydefault}{\mddefault}{\updefault}$t$}}}
\end{picture}

%% file: Versal8.pstex_t
\begin{picture}(0,0)%
\includegraphics{Versal8.pstex}%
\end{picture}%
\setlength{\unitlength}{3947sp}%
\begingroup\makeatletter\ifx\SetFigFont\undefined%
\gdef\SetFigFont#1#2#3#4#5{%
  \reset@font\fontsize{#1}{#2pt}%
  \fontfamily{#3}\fontseries{#4}\fontshape{#5}%
  \selectfont}%
\fi\endgroup%
\begin{picture}(5953,1306)(601,-3710)
\put(601,-2611){\makebox(0,0)[lb]{\smash{\SetFigFont{12}{14.4}{\familydefault}{\mddefault}{\updefault}$y$}}}
\put(6301,-2536){\makebox(0,0)[lb]{\smash{\SetFigFont{12}{14.4}{\familydefault}{\mddefault}{\updefault}$x$}}}
\put(676,-3661){\makebox(0,0)[lb]{\smash{\SetFigFont{12}{14.4}{\familydefault}{\mddefault}{\updefault}$t$}}}
\put(2326,-3661){\makebox(0,0)[lb]{\smash{\SetFigFont{12}{14.4}{\familydefault}{\mddefault}{\updefault}$t'$}}}
\put(4126,-3661){\makebox(0,0)[lb]{\smash{\SetFigFont{12}{14.4}{\familydefault}{\mddefault}{\updefault}$t''$}}}
\end{picture}

%% file: Versal9.pstex_t
\begin{picture}(0,0)%
\includegraphics{Versal9.pstex}%
\end{picture}%
\setlength{\unitlength}{3947sp}%
\begingroup\makeatletter\ifx\SetFigFont\undefined%
\gdef\SetFigFont#1#2#3#4#5{%
  \reset@font\fontsize{#1}{#2pt}%
  \fontfamily{#3}\fontseries{#4}\fontshape{#5}%
  \selectfont}%
\fi\endgroup%
\begin{picture}(5440,2684)(1039,-3014)
\put(1276,-1336){\makebox(0,0)[lb]{\smash{\SetFigFont{12}{14.4}{\familydefault}{\mddefault}{\updefault}$\Phi_1$}}}
\put(3901,-1486){\makebox(0,0)[lb]{\smash{\SetFigFont{12}{14.4}{\familydefault}{\mddefault}{\updefault}$\Phi_2$}}}
\put(4951,-1486){\makebox(0,0)[lb]{\smash{\SetFigFont{12}{14.4}{\familydefault}{\mddefault}{\updefault}$\Phi_3$}}}
\end{picture}